\documentclass[letterpaper,10pt,conference]{ieeeconf}  

\IEEEoverridecommandlockouts                              
\usepackage{amsmath}
\usepackage{amssymb}
\usepackage{algorithm}  
\usepackage{algorithmicx}  
\usepackage{algpseudocode}  
\usepackage{graphicx}
\usepackage{float}
\usepackage{frenchineq}

\usepackage{subcaption}
\usepackage{booktabs}

\usepackage{hyperref}
\usepackage{cleveref}
\usepackage[colorinlistoftodos,bordercolor=orange,backgroundcolor=orange!20,linecolor=orange,textsize=scriptsize]{todonotes}

\newcommand{\R}{\mathbb{R}}

\newcommand{\U}{\mathcal{U}}
\newcommand{\Y}{\mathcal{Y}}

\newcommand{\Lc}{\mathcal{L}}

\newcommand{\ov}[1]{\overline{#1}}
\newcommand{\norm}[1]{\left\lVert #1 \right\rVert}
\newcommand{\abs}[1]{\left\lvert #1 \right\rvert}

\makeatletter
\let\NAT@parse\undefined
\makeatother
\usepackage[numbers]{natbib}

\newtheorem{proposition}{Proposition}[section]

\newtheorem{assumption}{Assumption}[section]
\newtheorem{remark}{Remark}[section]

\title{\LARGE \bf PINNs in PDE Constrained Optimal Control Problems: Direct vs Indirect Methods}

\author{  Zhen Zhang, Shanqing Liu, Alessandro Alla, J\'er\^ome Darbon, George Em Karniadakis
	\thanks{Zhen Zhang, Shanqing Liu, J\'er\^ome Darbon, and George Em Karniadakis are with Division of Applied Mathematics, Brown University, {\tt\footnotesize \{zhen\_zhang2, shanqing\_liu, jerome\_darbon, george\_karniadakis\}@brown.edu}}%
    \thanks{Alessandro Alla is with Department of Mathematics at Sapienza, Università di Roma, {\tt\footnotesize alessandro.alla@uniroma1.it} }
        }

\begin{document}
	
	\maketitle
	\thispagestyle{empty} 
	\pagestyle{empty}
	\begin{abstract} 
    We study physics-informed neural networks (PINNs) as numerical tools for the optimal control of semilinear partial differential equations. We first recall the classical direct and indirect viewpoints for optimal control of PDEs, and then present two PINN formulations: a direct formulation based on minimizing the objective under the state constraint, and an indirect formulation based on the first-order optimality system. 
    For a class of semilinear parabolic equations, we derive the state equation, the adjoint equation, and the stationarity condition in a form consistent with continuous-time Pontryagin-type optimality conditions. We then specialize the framework to an Allen-Cahn control problem and compare three numerical approaches: (i) a discretize-then-optimize adjoint method, (ii) a direct PINN, and (iii) an indirect PINN. 
    Numerical results show that the PINN parameterization has an implicit regularizing effect, in the sense that it tends to produce smoother control profiles. They also indicate that the indirect PINN more faithfully preserves the PDE contraint and optimality structure and yields a more accurate neural approximation than the direct PINN.
\end{abstract}

	\section{INTRODUCTION}

Optimal control problems constrained by partial differential equations (PDEs) arise in a broad range of applications, including heat transfer, diffusion-reaction processes, fluid flows, phase-field dynamics, and distributed parameter systems more generally \cite{hinze2009, hinze2011discretization}. In these problems, the control acts on an infinite-dimensional dynamical system, while the objective functional must be optimized under the constraint that the state satisfies the governing PDE. This intrinsic coupling between optimization and dynamics leads to substantial computational and analytical difficulties, especially for nonlinear equations, fine space-time discretizations, control constraints, and settings in which repeated forward and adjoint solves are required \cite{troltzsch2010,hinze2009}. For time-dependent systems, the difficulty is compounded by the need to propagate information both forward through the state equation and backward through the adjoint equation, often on large computational grids.

A standard viewpoint distinguishes between \emph{direct} and \emph{indirect} approaches~\cite{cohen2015approximation, zhang2024sequential}. Direct methods discretize the state equation and then minimize the resulting finite-dimensional objective under the discrete dynamics. On the other hand, indirect methods derive, initially, the first-order necessary optimality conditions, typically a coupled state adjoint system together with a stationary condition, and then solve this optimality system numerically. Both viewpoints are classical in PDE-constrained optimization, and each has complementary advantages: direct methods are often more robust at the optimization level, whereas indirect methods more explicitly exploit the variational structure of the control problem.

In recent years, physics-informed neural networks (PINNs) have emerged as a flexible framework for the approximation of PDE solutions and related inverse problems \cite{raissi2019,karniadakis2021}. Their appeal in optimal control is natural, and is fruitful in the literature (see e.g~\cite{garcia2023control,mowlavi2023,barry2025,na2024physics,barry2025physics,meng2022sympocnet,zhang2025time}). The state, control, and possibly the adjoint can be represented by neural networks, while the governing equations, boundary conditions, and optimality relations are enforced through residual losses. This makes PINNs particularly attractive when one seeks low-dimensional parametric representations of the control and state, or when mesh-free training over the space-time domain is desirable. At the same time, for control problems it remains important to understand how different PINN formulations relate to the underlying optimality system, and whether this distinction has a tangible effect on the quality of the computed control.
We also note that PINNs have been used in the context of control problems with unknown state equations in \cite{alla2025pinn}. In that work, the contribution of PINNs lies in first identifying the underlying state equation and subsequently using it for control within an indirect approach.

In this paper, we focus on 
semilinear parabolic optimal control problems and study two PINN formulations for their numerical solution. The first is a \emph{direct} formulation, in which the state and control networks are trained by minimizing a loss combining the PDE residual, the boundary and initial conditions, and the control objective. The second is an \emph{indirect} formulation, in which the training loss is built from the first-order optimality system, 
that is a system of equations consist of first-order system of  the state equation and the adjoint equation, and the first-order optimality condition with respect to the control.
We compare these two neural formulations with a classical discretize-then-optimize adjoint method on a one-dimensional Allen-Cahn 
control problem.

This paper is organized as follows.
First, we present a concise control-theoretic formulation of the semilinear parabolic optimal control problem and of its associated first-order optimality system. Second, we show how this structure leads to two distinct PINN formulations, namely a direct objective-based approach and an indirect KKT/Pontryagin-type approach. Third, through a numerical comparison on the Allen-Cahn equation, we show that the indirect formulation yields a more informative training signal and produces the most accurate neural approximation among the methods considered. These results suggest that, for PDE-constrained optimal control, enforcing the optimality system itself may be more effective than relying solely on objective minimization within the PINN framework. 
From a practical perspective, the results indicate that successful PINN training benefits strongly from second-order optimization and double-precision arithmetic. They also show that PINNs offer a user-friendly framework and can serve as convenient initializers for other numerical approaches, including adjoint and indirect PINNs.

    \section{Preliminary} 

    \subsection{PDE-constrained optimal control} 
    Let $\Omega\subset\R^d$ be a bounded  domain, let $T>0$, and define $Q:=\Omega\times(0,T), 
\Sigma:=\partial\Omega\times(0,T)$. 
We consider the semilinear parabolic state equation
\begin{equation}
\label{eq:abstract_state}
\begin{cases}
\partial_t y - \nu\Delta y + f(y) = \mathcal{B}u, & \text{in } Q,\\
\partial_n y = 0, & \text{on } \Sigma,\\
y(\cdot,0)=y_0, & \text{in } \Omega,
\end{cases}
\end{equation}
where $\nu>0$ is the diffusion coefficient, $u$ is the distributed control, $\mathcal{B}$ is a bounded control operator, and $f$ is a possibly nonlinear reaction term.

A standard control space is $\U := L^2(Q)$, and the natural weak state space is
\[
\Y := L^2(0,T;H^1(\Omega))\cap H^1(0,T;H^1(\Omega)^*).
\]
The optimal control problem consists in minimizing
\begin{equation}
\label{eq:abstract_cost}
J(y,u) := \frac{\beta_T}{2}\norm{y(\cdot,T)-y_d}_{L^2(\Omega)}^2
+ \frac{\beta_Q}{2}\norm{u}_{L^2(Q)}^2,
\end{equation}
subject to \eqref{eq:abstract_state}, where $y_d\in L^2(\Omega)$ is a desired terminal target and $\beta_T,\beta_Q>0$ are weights.  

We shall use the following standard structural assumption.

\begin{assumption}
\label{ass:nonlinearity}
The nonlinearity $f:\R\to\R$ is of class $C^1$, locally Lipschitz, and its derivative is bounded from below, i.e., there exists $c_f\ge 0$ such that
\[
f'(s)\ge -c_f, \qquad \forall s\in\R.
\]
The initial datum satisfies $y_0\in L^2(\Omega)$, and $\mathcal{B}:L^2(Q)\to L^2(Q)$ is bounded.
\end{assumption}

\begin{proposition}
\label{prop:state_wellposed}
Under Assumption~\ref{ass:nonlinearity}, for every $u\in\U$, the state equation \eqref{eq:abstract_state} admits a unique weak solution $y=S(u)\in\Y$. Moreover, the control-to-state map $S:\U\to\Y$ is well defined and locally Fr\'echet differentiable.
\end{proposition} 
\begin{remark}
    This class of control problems, with the state equation given by \eqref{eq:abstract_state}, includes, for instance, the Allen-Cahn equation $f(y) = y^3-y$, the degenerate Zeldovich equation $f(y) = y^3-y^2$, and the Zeldovich-Frank-Kamenetskii $f(y) = y(y-1)(y-a) $ with $a\in (0,1)\subset\mathbb{R}$. Numerical studies of controlled problems for these equationscan be found in the literature, see e.g. \cite{altmuller10,alla23}.
\end{remark}
\begin{remark}
For the semilinear parabolic problems considered here, Proposition~\ref{prop:state_wellposed} follows from standard monotonicity and energy estimates. The paper does not rely on a new well-posedness result, but rather on the corresponding optimality system and its numerical approximation.
\end{remark}

\subsection{Direct and indirect numerical viewpoints} 

Up to rare cases, problem of the form~\eqref{eq:abstract_cost} can not be solved exactly. In the numerical approximation aspects, two classical viewpoints coexist in optimal control.

\emph{Direct methods} optimize the reduced functional
\[
j(u):=J(S(u),u),
\]
possibly after discretizing the state equation first. In the PDE literature, this includes both optimize-then-discretize and discretize-then-optimize workflows, depending on whether one derives the first-order conditions before or after numerical discretization.

\emph{Indirect methods}, instead, first derive the first-order optimality system, typically through a Lagrange multiplier or Pontryagin maximum principle argument, and then solve the coupled state-adjoint-stationarity system. Indirect methods are often more structured and informative, but they require access to the adjoint dynamics and may be more delicate numerically.

For the present paper, this distinction is especially relevant because it leads to two different PINN formulations. The direct PINN minimizes an objective-based (that is the cost functional of the control problem,  e.g. \eqref{eq:abstract_cost}) loss under the state constraint, whereas the indirect PINN enforces the first-order optimality system itself.

\subsection{Physics-informed neural networks} 

Physics-informed neural networks (PINNs) approximate unknown functions by neural networks and train them by minimizing residuals of governing equations, boundary conditions, and data or objective terms. 
For PDE-constrained optimization, one may represent the state, control, and possibly the adjoint by neural networks,
\[
Y_\theta(x,t)\approx y(x,t),
\ 
U_\psi(x,t)\approx u(x,t),
\
\Lambda_\phi(x,t)\approx \lambda(x,t) \ ,
\]
and define a composite loss from automatic-differentiation-based residuals.  

The attraction of this approach is twofold. First, the formulation is continuous in space and time, which avoids committing a priori to a specific grid in the optimization variables. Second, a system of  equations  derived from the maximum principle can be coupled in a unified loss, which makes PINNs particularly appealing for indirect formulations involving both state and adjoint equations.

At the same time, the accuracy of PINNs depends on network expressivity, loss balancing, optimization, and collocation design. In particular, for optimal control problems, the choice between direct and indirect formulations can significantly affect optimization behavior and final solution quality.

\section{PINN FORMULATIONS FOR OPTIMAL CONTROL OF SEMILINEAR PDES} 

\subsection{Problem formulation}

In this paper, we focus on the unconstrained distributed control problem
\begin{equation}
\label{eq:ocp_general}
\min_{(y,u)} J(y,u)
\quad
\text{subject to } (y,u) \text{ satisfying } \eqref{eq:abstract_state}.
\end{equation}
An admissible pair is any $(y,u)\in\Y\times\U$ satisfying the state equation in the weak sense. Since the control cost in \eqref{eq:abstract_cost} is coercive in $u$, standard compactness arguments yield existence of at least one optimal pair for the class of problems considered here (see e.g \cite{troltzsch2010}). We do not assume global convexity of the reduced problem because the semilinearity may render the optimization nonconvex. 

To derive first-order necessary conditions, introduce the Lagrangian
\begin{equation}
\label{eq:lagrangian_general}
\begin{aligned}
\Lc(y,u,\lambda)
&:=  \frac{\beta_T}{2}\norm{y(\cdot,T)-y_d}_{L^2(\Omega)}^2
+ \frac{\beta_Q}{2}\norm{u}_{L^2(Q)}^2 \\
&+ \int_0^T\!\int_\Omega \lambda\,\bigl(\partial_t y-\nu\Delta y+f(y)-\mathcal{B}u\bigr)\,dx\,dt.
\end{aligned}
\end{equation}
Here $\lambda$ is the adjoint variable associated with the state constraint.

\subsection{First-order optimality system and maximum principle}

Assume that $(\ov y,\ov u)$ is a locally optimal pair with sufficient regularity, and let $\ov\lambda$ be the associated adjoint. Taking first variations in \eqref{eq:lagrangian_general} yields the following continuous first-order system:
\begin{equation}
\label{eq:state_opt}
\begin{cases}
\partial_t \ov y - \nu\Delta \ov y + f(\ov y)=\mathcal{B}\ov u, & \text{in } Q,\\
\partial_n \ov y = 0, & \text{on } \Sigma,\\
\ov y(\cdot,0)=y_0, & \text{in } \Omega,
\end{cases}
\end{equation}
\begin{equation}
\label{eq:adjoint_opt}
\begin{cases}
-\partial_t \ov\lambda - \nu\Delta \ov\lambda + f'(\ov y)\,\ov\lambda = 0, & \text{in } Q,\\
\partial_n \ov\lambda = 0, & \text{on } \Sigma,\\
\ov\lambda(\cdot,T) = -\beta_T\bigl(\ov y(\cdot,T)-y_d\bigr), & \text{in } \Omega,
\end{cases}
\end{equation}
and the first-order optimality condition with respect to the control
\begin{equation}
\label{eq:stationarity_general}
\beta_Q\,\ov u - \mathcal{B}^*\ov\lambda = 0
\qquad \text{in } Q.
\end{equation}
Equivalently,
\begin{equation}
\label{eq:feedback_general}
\ov u = \beta_Q^{-1}\mathcal{B}^*\ov\lambda.
\end{equation} 

\begin{proposition}
\label{prop:foc}
Suppose that Assumption~\ref{ass:nonlinearity} holds and that $(\ov y,\ov u)$ is a local minimizer of \eqref{eq:ocp_general}. Then, there exists an adjoint variable $\ov\lambda$ such that \eqref{eq:state_opt}--\eqref{eq:stationarity_general} hold.
\end{proposition}


\subsection{Direct PINN formulation} 

In the direct formulation, the state and control are parameterized by neural networks,
\[
y\approx Y_\theta, \qquad u\approx U_\psi,
\]
and the loss is built from the state residual together with the objective. Let
\[
r_y(x,t;\theta,\psi)
:= \partial_t Y_\theta - \nu\Delta Y_\theta + f(Y_\theta) - \mathcal{B}U_\psi.
\] 
A generic direct PINN loss takes the form
\begin{equation}
\label{eq:loss_direct_general}
\begin{aligned}
\mathcal{L}_{\mathrm{dir}}(\theta,\psi)
:={}& w_{\mathrm{res}}\,\mathcal{L}_{\mathrm{res}}^y
+ w_{\mathrm{bc}}\,\mathcal{L}_{\mathrm{bc}}^y
+ w_{\mathrm{ic}}\,\mathcal{L}_{\mathrm{ic}}^y \\
&+ \frac{\beta_T}{2}\,\mathcal{L}_{T}
+ \frac{\beta_Q}{2}\,\mathcal{L}_{u},
\end{aligned}
\end{equation}
where each term is estimated from collocation points. For example, if
$\{(x_j,t_j)\}_{j=1}^{N_{\mathrm{int}}}\subset Q$ denotes interior collocation
points, $\{t_j^{\mathrm{bc}}\}_{j=1}^{N_{\mathrm{bc}}}\subset(0,T)$ boundary collocation
points, and $\{x_i\}_{i=1}^{N_T}\subset\Omega$ terminal samples, then
\begin{align*}
\mathcal{L}_{\mathrm{res}}^y
&\approx \frac{1}{N_{\mathrm{int}}}\sum_{j=1}^{N_{\mathrm{int}}}
\abs{r_y(x_j,t_j;\theta,\psi)}^2,\\
\mathcal{L}_{\mathrm{bc}}^y
&\approx \frac{1}{N_{\mathrm{bc}}}\sum_{j=1}^{N_{\mathrm{bc}}}
\abs{\partial_n Y_\theta(\cdot,t_j^{\mathrm{bc}})}^2,\\
\mathcal{L}_{\mathrm{ic}}^y
&\approx \frac{1}{N_0}\sum_{i=1}^{N_0}
\abs{Y_\theta(x_i,0)-y_0(x_i)}^2,\\
\mathcal{L}_{T}
&\approx \frac{|\Omega|}{N_T}\sum_{i=1}^{N_T}
\abs{Y_\theta(x_i,T)-y_d(x_i)}^2,\\
\mathcal{L}_{u}
&\approx \frac{|Q|}{N_{\mathrm{int}}}\sum_{j=1}^{N_{\mathrm{int}}}
\abs{U_\psi(x_j,t_j)}^2.
\end{align*}
Thus $\mathcal{L}_{\mathrm{bc}}^y$ and $\mathcal{L}_{\mathrm{ic}}^y$ enforce the boundary and initial conditions, while $\mathcal{L}_T$ and $\mathcal{L}_u$ approximate the two terms of the cost functional.


The direct formulation is conceptually simple and does not require the adjoint variable explicitly. However, the objective appears only through scalar integral terms, and this may provide a weaker optimization signal than enforcing the full optimality system.

\subsection{Indirect PINN formulation}

In the indirect formulation, the state, control, and adjoint are trained so as to satisfy the first-order system \eqref{eq:state_opt}--\eqref{eq:stationarity_general}. Introduce network approximations
\[
y\approx Y_\theta,
\qquad
u\approx U_\psi,
\qquad
\lambda\approx \Lambda_\phi,
\]
and define the residuals
\begin{align}
\label{eq:residual_state_general}
r_y(x,t;\theta,\psi)
&:= \partial_t Y_\theta - \nu\Delta Y_\theta + f(Y_\theta) - \mathcal{B}U_\psi,\\
\label{eq:residual_adj_general}
r_\lambda(x,t;\theta,\phi)
&:= -\partial_t \Lambda_\phi - \nu\Delta \Lambda_\phi + f'(Y_\theta)\Lambda_\phi,\\
\label{eq:residual_st_general}
r_{\mathrm{st}}(x,t;\psi,\phi)
&:= \beta_Q U_\psi - \mathcal{B}^*\Lambda_\phi.
\end{align}
A generic indirect PINN loss is then
\begin{equation}
\label{eq:loss_indirect_general}
\begin{aligned}
\mathcal{L}_{\mathrm{ind}}(\theta,&\psi,\phi)
:= w_y\,\mathcal{L}_{\mathrm{res}}^y
+ w_\lambda\,\mathcal{L}_{\mathrm{res}}^\lambda
+ w_{\mathrm{st}}\,\mathcal{L}_{\mathrm{st}} \\
&+ w_{\mathrm{bc}}^y\,\mathcal{L}_{\mathrm{bc}}^y
+ w_{\mathrm{bc}}^\lambda\,\mathcal{L}_{\mathrm{bc}}^\lambda
+ w_{\mathrm{ic}}\,\mathcal{L}_{\mathrm{ic}}^y
+ w_T\,\mathcal{L}_{T}^{\lambda}.
\end{aligned}
\end{equation}
Using the same sets of collocation points as above, one may write
\begin{align*}
\mathcal{L}_{\mathrm{res}}^\lambda
&\approx \frac{1}{N_{\mathrm{int}}}\sum_{j=1}^{N_{\mathrm{int}}}
\abs{r_\lambda(x_j,t_j;\theta,\phi)}^2,\\
\mathcal{L}_{\mathrm{st}}
&\approx \frac{1}{N_{\mathrm{int}}}\sum_{j=1}^{N_{\mathrm{int}}}
\abs{r_{\mathrm{st}}(x_j,t_j;\psi,\phi)}^2,\\
\mathcal{L}_{\mathrm{bc}}^\lambda
&\approx \frac{1}{N_{\mathrm{bc}}}\sum_{j=1}^{N_{\mathrm{bc}}}
\abs{\partial_n \Lambda_\phi(\cdot,t_j^{\mathrm{bc}})}^2,\\
\mathcal{L}_{T}^{\lambda}
&\approx \frac{|\Omega|}{N_T}\sum_{i=1}^{N_T}
\abs{\Lambda_\phi(x_i,T)+\beta_T\bigl(Y_\theta(x_i,T)-y_d(x_i)\bigr)}^2.
\end{align*}
Here $\mathcal{L}_{\mathrm{res}}^\lambda$ penalizes the adjoint PDE residual, $\mathcal{L}_{\mathrm{st}}$ penalizes the stationarity residual, and $\mathcal{L}_{T}^{\lambda}$ enforces the terminal condition for the adjoint.

The indirect formulation more closely mirrors the structure of classical indirect methods. In particular, a vanishing indirect PINN loss implies that the neural approximations satisfy the continuous first-order system at the collocation points. This typically yields a stronger and more informative training signal than the direct formulation, at the price of a larger and more tightly coupled loss.

    \section{Numerical results}

\subsection{Problem formulation}

We consider the Allen--Cahn equation on $\Omega = [0,1]$ with homogeneous Neumann boundary conditions:
\begin{equation}\label{eq:ac}
  y_t = \varepsilon^2\, y_{xx} - (y^3 - y) + u(x,t),
  \qquad (x,t)\in\Omega\times(0,T],
\end{equation}
subject to
\begin{alignat}{2}
  y_x(0,t) &= y_x(1,t) = 0, &\qquad& t\in(0,T], \label{eq:ac_bc}\\
  y(x,0)   &= \cos(\pi x),   &\qquad& x\in\Omega, \label{eq:ac_ic}
\end{alignat}
where $\varepsilon = 0.01$ is the interface thickness parameter, $T = 3.0$,
and $u(x,t)$ is the distributed control.

The optimal control problem seeks to minimize the cost functional
\begin{equation}\label{eq:ac_obj}
  J(y,u) = \frac{\beta_T}{2}\int_0^1 y(x,T)^2\,dx
          + \frac{\beta_Q}{2}\int_0^T\!\!\int_0^1 u(x,t)^2\,dx\,dt,
\end{equation}
subject to the state equation~\eqref{eq:ac}--\eqref{eq:ac_ic},
with $\beta_T = 1$ and $\beta_Q = 10^{-3}$. In the form of~\eqref{eq:abstract_state}, this corresponds to $\nu = \varepsilon^2$, $\mathcal{B}=I$,
$f(y)=y^3-y$, and $y_d\equiv 0$. The terminal cost drives the state toward
$y(x,T)=0$, while the control penalty regularizes the magnitude of $u$.

For all three methods, we use the same reference discretization of the state
PDE when reporting trajectories and objective values. 
The spatial interval is
partitioned by the uniform grid $x_i=(i-1)\Delta x$, $i=1,\dots,N$, with
$N=513$ and $\Delta x=1/512$, and the time interval is discretized by
$t_n=n\Delta t$, $n=0,\dots,N_t$, with $\Delta t=0.05$ and $N_t=60$. 
The 
second derivative is approximated by the standard centered finite-difference
operator $D_{xx}$ together with Neumann ghost-point closures,
$Y_0^n=Y_2^n$ and $Y_{N+1}^n=Y_{N-1}^n$, so that
$\partial_x y(0,t)=\partial_x y(1,t)=0$ is enforced to second order. 
The 
resulting semi-discrete system is advanced by the classical explicit four-stage
Runge-Kutta (RK4) method. The same grid and quadrature are used to evaluate all
controls produced by the PINN and adjoint approaches, ensuring that the final
comparison is performed on a common discrete baseline.

\subsection{Method 1: Adjoint method (discretize-then-optimize)}

The adjoint baseline uses exactly the discretization described above.  
The 
control is represented by a piecewise-constant-in-time array
$\mathbf{U}=\bigl(U_i^n\bigr)_{1\le i\le N,\,0\le n\le N_t-1}$, with
$U_i^n\approx u(x_i,t_n)$, yielding $513\times 60=30{,}780$ optimization
degrees of freedom. Let $\mathbf{Y}=\bigl(Y_i^n\bigr)$ denote the associated
discrete state produced by the RK4 solver.
The discrete adjoint of the RK4 time-stepper is derived analytically and
used to compute the exact gradient $\nabla_{\mathbf{U}}J$ in a single backward sweep.
Optimization is performed with a second-order optimizer SSBroyden proposed in \cite{optimizer}. The cost functional~\eqref{eq:ac_obj} is approximated by 
\begin{align}
  J_T^h &= \frac{\beta_T}{2}\,\Delta x \sum_{i=1}^{N} \bigl(Y_i^{N_t}\bigr)^2, \\
  J_Q^h &= \frac{\beta_Q}{2}\,\Delta x\,\Delta t \sum_{n=0}^{N_t-1}\sum_{i=1}^{N} \bigl(U_i^n\bigr)^2,
\end{align}
so that $J^h = J_T^h + J_Q^h$.

\subsection{Method 2: PINN -- direct formulation}
\label{sec:pinn_direct}
\paragraph{Neural network architecture}
Two independent multi-layer perceptrons (MLPs) are used:
\begin{itemize}
  \item $Y_\theta(x,t)$:\ state network, mapping $(x,t)\mapsto y$,
  \item $U_\psi(x,t)$:\ control network, mapping $(x,t)\mapsto u$.
\end{itemize}
Each MLP has input dimension~$2$, output dimension~$1$, two hidden layers of
width~$32$, and $\tanh$ activation, yielding $1{,}185$ trainable parameters
per network ($2{,}370$ total).
The weights are initialized with Xavier uniform initialization.
\paragraph{Loss function}
The PINN is trained by minimizing the composite loss
\begin{equation}\label{eq:pinn_direct_loss}
\begin{aligned}
  \mathcal{L}_{\mathrm{dir}}(\theta,\psi)
    ={}& w_{\mathrm{res}}\,\mathcal{L}_{\mathrm{res}}^y
    + w_{\mathrm{bc}}\,\mathcal{L}_{\mathrm{bc}}^y
    + w_{\mathrm{ic}}\,\mathcal{L}_{\mathrm{ic}}^y \\
    &+ \frac{\beta_T}{2}\,\mathcal{L}_T
    + \frac{\beta_Q}{2}\,\mathcal{L}_u,
\end{aligned}
\end{equation}
where, for the Allen--Cahn equation,
\[
r_y(x,t;\theta,\psi)
:= \partial_t Y_\theta - \varepsilon^2\partial_{xx}Y_\theta
   + \bigl(Y_\theta^3-Y_\theta\bigr) - U_\psi.
\]
The loss components are
\begin{align}
  \mathcal{L}_{\mathrm{res}}^y &= \frac{1}{N_{\mathrm{int}}} \sum_{j=1}^{N_{\mathrm{int}}}
    \bigl| r_y(x_j,t_j;\theta,\psi) \bigr|^2, \\[4pt]
  \mathcal{L}_{\mathrm{bc}}^y &= \frac{1}{N_{\mathrm{bc}}} \sum_{j=1}^{N_{\mathrm{bc}}}
    \bigl( |\partial_xY_\theta(0,t_j)|^2 + |\partial_xY_\theta(1,t_j)|^2 \bigr), \\[4pt]
  \mathcal{L}_{\mathrm{ic}}^y &= \frac{1}{N} \sum_{i=1}^{N}
    \bigl| Y_\theta(x_i,0) - y_0(x_i) \bigr|^2, \\[4pt]
  \mathcal{L}_T &= \frac{1}{N} \sum_{i=1}^{N}
    \bigl| Y_\theta(x_i,T) \bigr|^2, \\[4pt]
  \mathcal{L}_u &= \frac{T}{N_{\mathrm{int}}} \sum_{j=1}^{N_{\mathrm{int}}}
    \bigl| U_\psi(x_j,t_j) \bigr|^2.
\end{align}
Thus $\mathcal{L}_T$ and $\mathcal{L}_u$ are Monte Carlo approximations of the
two terms in~\eqref{eq:ac_obj}, and, since $|\Omega|=1$, no additional spatial
measure factor appears. 
All collocation points are generated without labeled
data, and the training signal comes entirely from the PDE, boundary/initial
conditions, and cost terms.
The loss weights are set to
$w_{\mathrm{res}} = w_{\mathrm{bc}} = w_{\mathrm{ic}} = 1$.


\paragraph{Collocation points}
$N_{\mathrm{int}} = 20{,}000$ interior points are sampled uniformly in
$[0,1]\times(0,T)$, and $N_{\mathrm{bc}} = 512$ boundary times are sampled
uniformly in $(0,T)$.
The initial condition and terminal condition are each evaluated on $N=513$
uniformly spaced grid points.
The collocation points are drawn once with a fixed random seed and held constant
throughout training (no resampling between epochs).

\paragraph{Optimizer}
Training proceeds in two phases:
\begin{enumerate}
  \item \emph{Adam warm-up}: $1{,}000$ steps with an initial learning rate of
    $10^{-3}$ and a step-decay schedule (factor $0.3$ every $200$ steps).
  \item \emph{SSBroyden}: $40$ outer epochs, each running the SSBroyden for up to $200$ iterations.
    At the start of each epoch, the collocation batch is randomly regenerated.
\end{enumerate}
All computations use \texttt{float64} precision throughout.

\subsection{Method 3: PINN -- indirect (KKT) formulation}
\label{sec:pinn_kkt}

\paragraph{Neural network architecture}
The same two-network architecture as the direct formulation
(Section~\ref{sec:pinn_direct}) is used:
$Y_\theta$ for the state and $U_\psi$ for the control/adjoint.
Since $\mathcal{B}=I$, the first-order optimality
condition reduces to $\beta_Q u - \lambda = 0$. We therefore define the
adjoint by construction as
\[
\Lambda_\psi(x,t):=\beta_Q U_\psi(x,t),
\]
so that $u(x,t)=U_\psi(x,t)$ and $\lambda(x,t)=\Lambda_\psi(x,t)$. In this
way the optimality condition is enforced exactly and does not need a separate
loss term.

\paragraph{Loss function}
Following~\eqref{eq:loss_indirect_general}, the indirect PINN minimizes
\begin{equation}\label{eq:pinn_kkt_loss}
\begin{aligned}
  \mathcal{L}_{\mathrm{ind}}(\theta,\psi)
  ={}& w_y\,\mathcal{L}_{\mathrm{res}}^y
    + w_\lambda\,\mathcal{L}_{\mathrm{res}}^\lambda
    + w_{\mathrm{bc}}^y\,\mathcal{L}_{\mathrm{bc}}^y \\
    &+ w_{\mathrm{bc}}^\lambda\,\mathcal{L}_{\mathrm{bc}}^\lambda
    + w_{\mathrm{ic}}\,\mathcal{L}_{\mathrm{ic}}^y
    + w_T\,\mathcal{L}_{T}^{\lambda}.
\end{aligned}
\end{equation} 
Here $\mathcal{L}_{\mathrm{res}}^y$, $\mathcal{L}_{\mathrm{bc}}^y$, and
$\mathcal{L}_{\mathrm{ic}}^y$ are exactly the same as in the direct
formulation. Since $\lambda=\Lambda_\psi=\beta_Q U_\psi$ is enforced by
construction, the specialized indirect loss omits the stationarity term
$\mathcal{L}_{\mathrm{st}}$. The remaining indirect terms are
\begin{align}
  \mathcal{L}_{\mathrm{res}}^\lambda &= \frac{1}{N_{\mathrm{int}}}
    \sum_{j=1}^{N_{\mathrm{int}}}
    \bigl| {-}\partial_t\Lambda_\psi - \varepsilon^2\partial_{xx}\Lambda_\psi
    + (3Y_\theta^2-1)\Lambda_\psi \bigr|^2_{(x_j,t_j)}, \\[4pt]
  \mathcal{L}_{\mathrm{bc}}^\lambda &= \frac{1}{N_{\mathrm{bc}}}
    \sum_{j=1}^{N_{\mathrm{bc}}}
    \bigl(|\partial_x\Lambda_\psi(0,t_j)|^2 + |\partial_x\Lambda_\psi(1,t_j)|^2\bigr), \\[4pt]
  \mathcal{L}_{T}^{\lambda} &= \frac{1}{N}\sum_{i=1}^{N}
    \bigl|\Lambda_\psi(x_i,T) + \beta_T Y_\theta(x_i,T)\bigr|^2.
\end{align}
Because $y_d\equiv 0$ and $\beta_T=1$ here, the last term enforces the
adjoint terminal condition $\lambda(x,T)=-y(x,T)$. All loss weights are set to
$w_y = w_\lambda = w_{\mathrm{bc}}^y = w_{\mathrm{bc}}^\lambda = w_{\mathrm{ic}} = w_T = 1$.


\paragraph{Collocation points and optimizer}
Same as the direct formulation: $N_{\mathrm{int}} = 20{,}000$ interior points, $N_{\mathrm{bc}} = 512$ boundary times, $1{,}000$ Adam warm-up steps followed by $40$ SSBroyden epochs (up to $200$ iterations each), all in
\texttt{float64}.

\subsection{Solution evaluation}

For all three methods, the reported state trajectory $y(x,t)$ is obtained by
feeding the converged control $u(x,t)$ into the RK4 solver
(Section~\ref{sec:pinn_direct} uses the same $N=513$, $\Delta t = 0.05$
discretization as the adjoint method) rather than using the PINN's direct
neural-network prediction.
This ensures that all methods are compared on the same numerical baseline and
that the reported state satisfies the discrete state equation exactly for the
given control.
The objective values $J_T$, $J_u$, and $J = J_T + J_u$ are likewise computed
from the solver trajectory using the trapezoidal rule.

\subsection{Discussion of numerical results}
\label{sec:result_discussion}

\begin{figure}[h!]
    \centering

    \begin{subfigure}[b]{0.7\columnwidth}
        \centering
        \includegraphics[width=\textwidth]{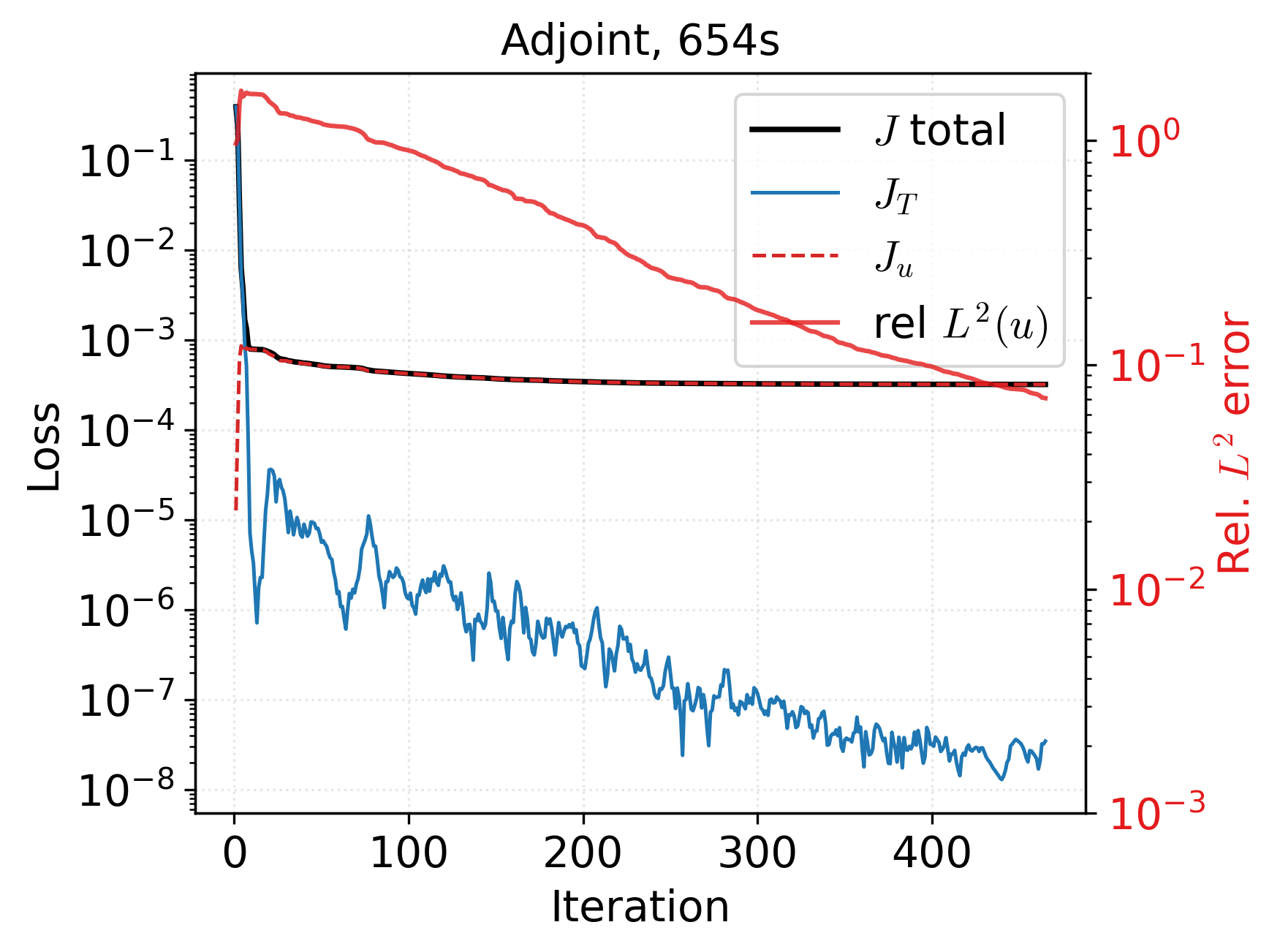}
        \label{fig:loss_adjoint}
    \end{subfigure}

    \vspace{0.1em}

    \begin{subfigure}[b]{0.7\columnwidth}
        \centering
        \includegraphics[width=\textwidth]{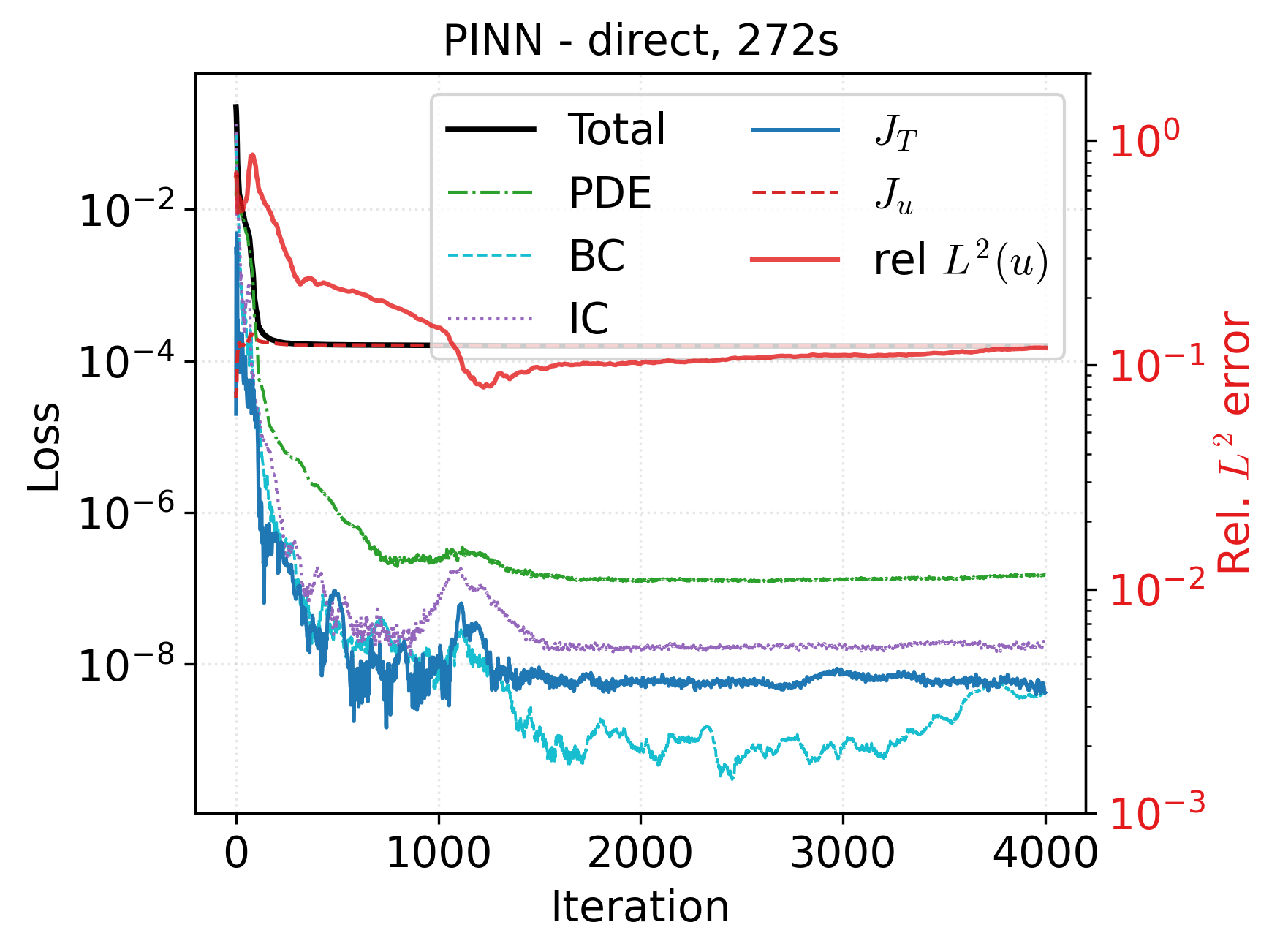}
        \label{fig:loss_pinn_direct}
    \end{subfigure}

    \vspace{0.1em}

    \begin{subfigure}[b]{0.7\columnwidth}
        \centering
        \includegraphics[width=\textwidth]{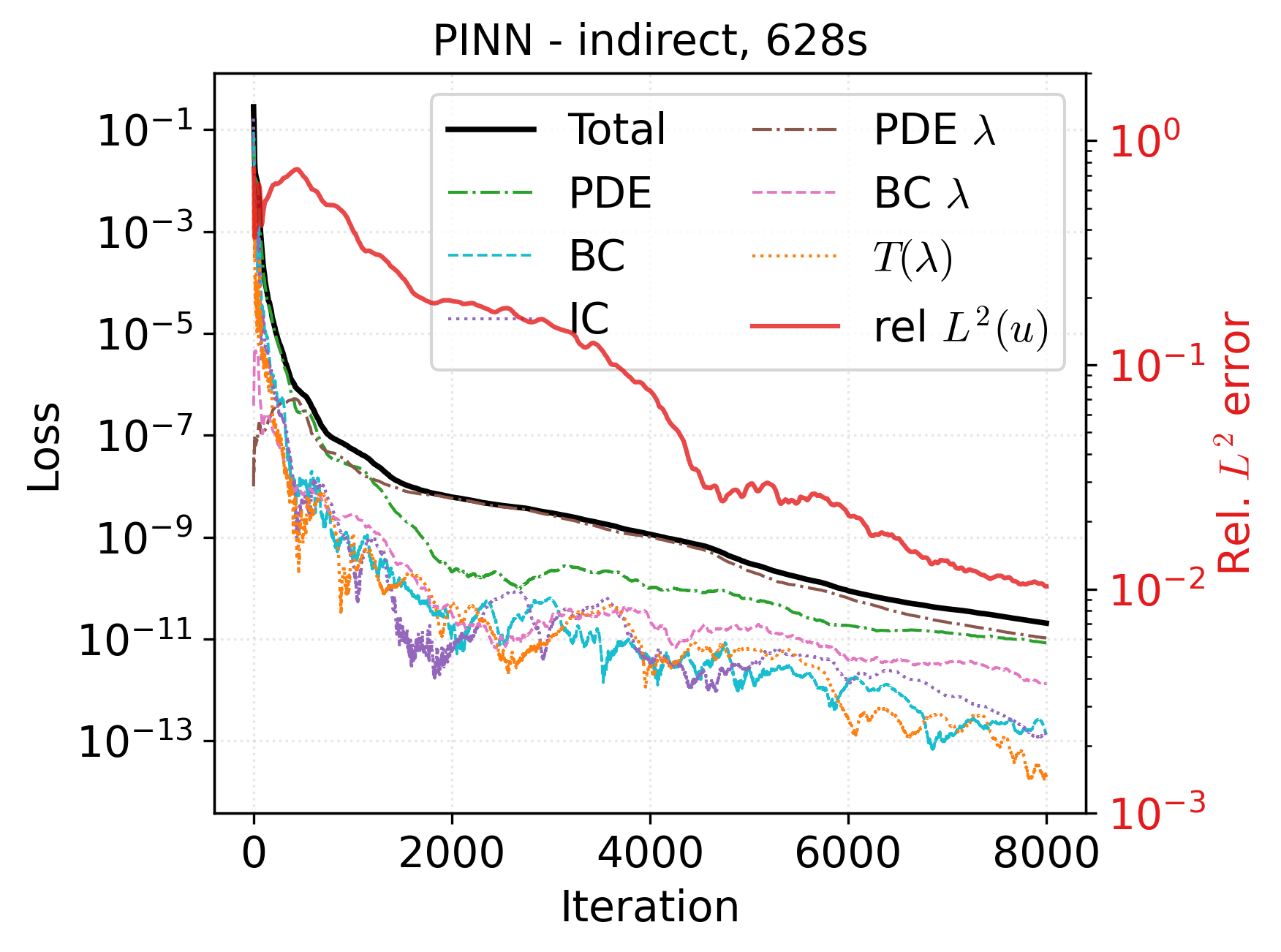}
        \label{fig:loss_pinn_indirect}
    \end{subfigure}
    \caption{Loss histories for the three approaches: adjoint optimization starting from scratch (top), direct PINN (middle), and indirect PINN (bottom). In the adjoint panel, the legend shows the total objective $\mathcal{J}$, the terminal tracking term $\mathcal{J}_T$, the control regularization term $\mathcal{J}_u$, and the relative $L_2$ error of the control, $\mathrm{rel}\,L^2(u)$. In the direct PINN panel, the legend additionally separates the PDE, boundary-condition (BC), and initial-condition (IC) residual losses, together with $\mathcal{J}_T$, $\mathcal{J}_u$, and $\mathrm{rel}\,L^2(u)$. In the indirect PINN panel, the legend further includes the residual losses associated with the adjoint equation ($\mathrm{PDE}\ \lambda$), adjoint boundary condition ($\mathrm{BC}\ \lambda$), and terminal optimality condition ($\mathcal{T}(\lambda)$), in addition to the state PDE/BC/IC residuals, the total loss, and $\mathrm{rel}\,L^2(u)$. For the PINN cases, only the SSBroyden phase is shown after the short Adam warm-up, and the total wall-clock time of each method is indicated in the panel title.}
    \label{fig:losses}
\end{figure}

\begin{figure}[!ht]
    \centering

    \begin{subfigure}[b]{0.7\columnwidth}
        \centering
        \includegraphics[width=\textwidth]{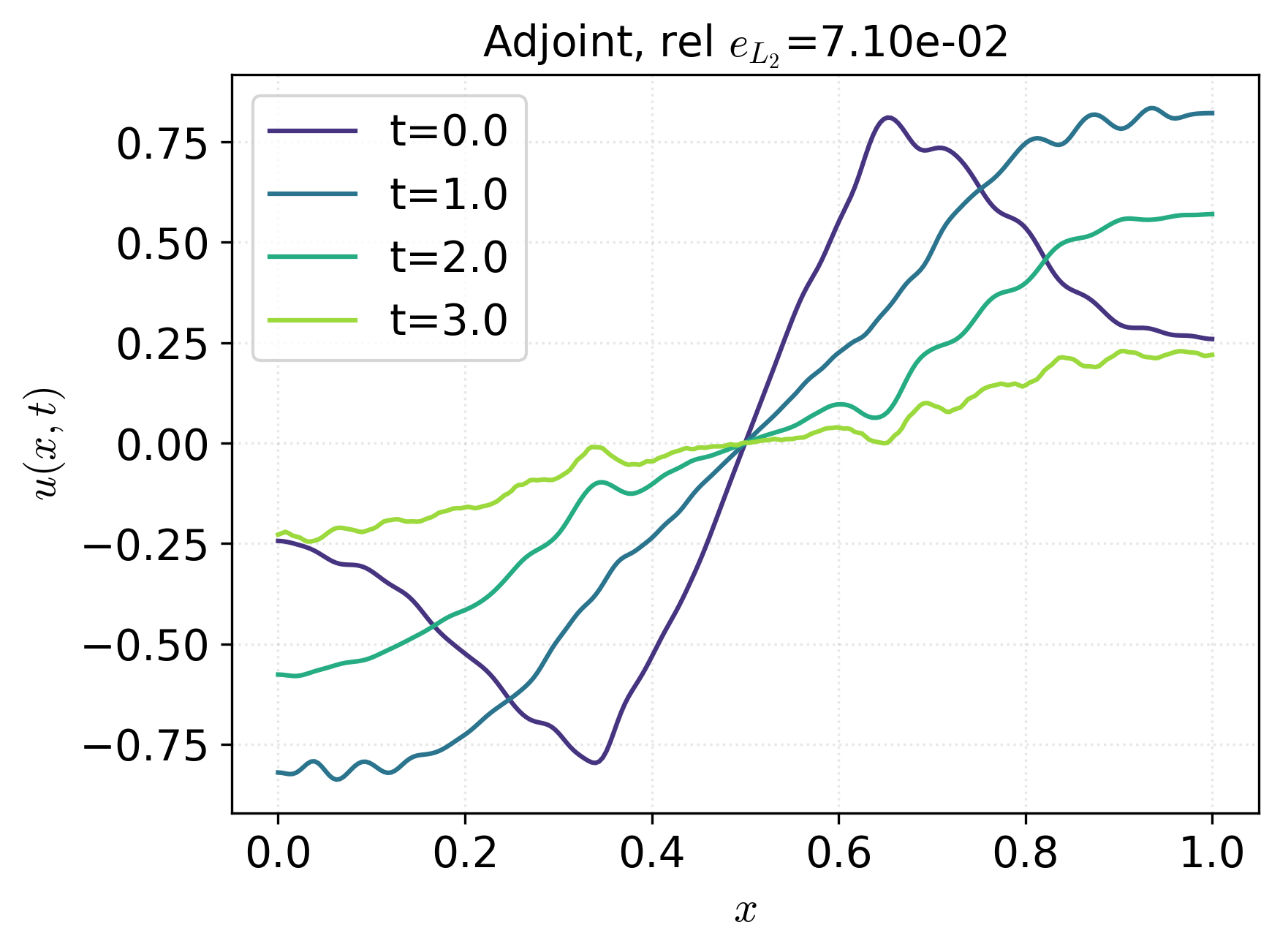}
        \label{fig:ctrl_adjoint}
    \end{subfigure}

    \vspace{0.1em}

    \begin{subfigure}[b]{0.7\columnwidth}
        \centering
        \includegraphics[width=\textwidth]{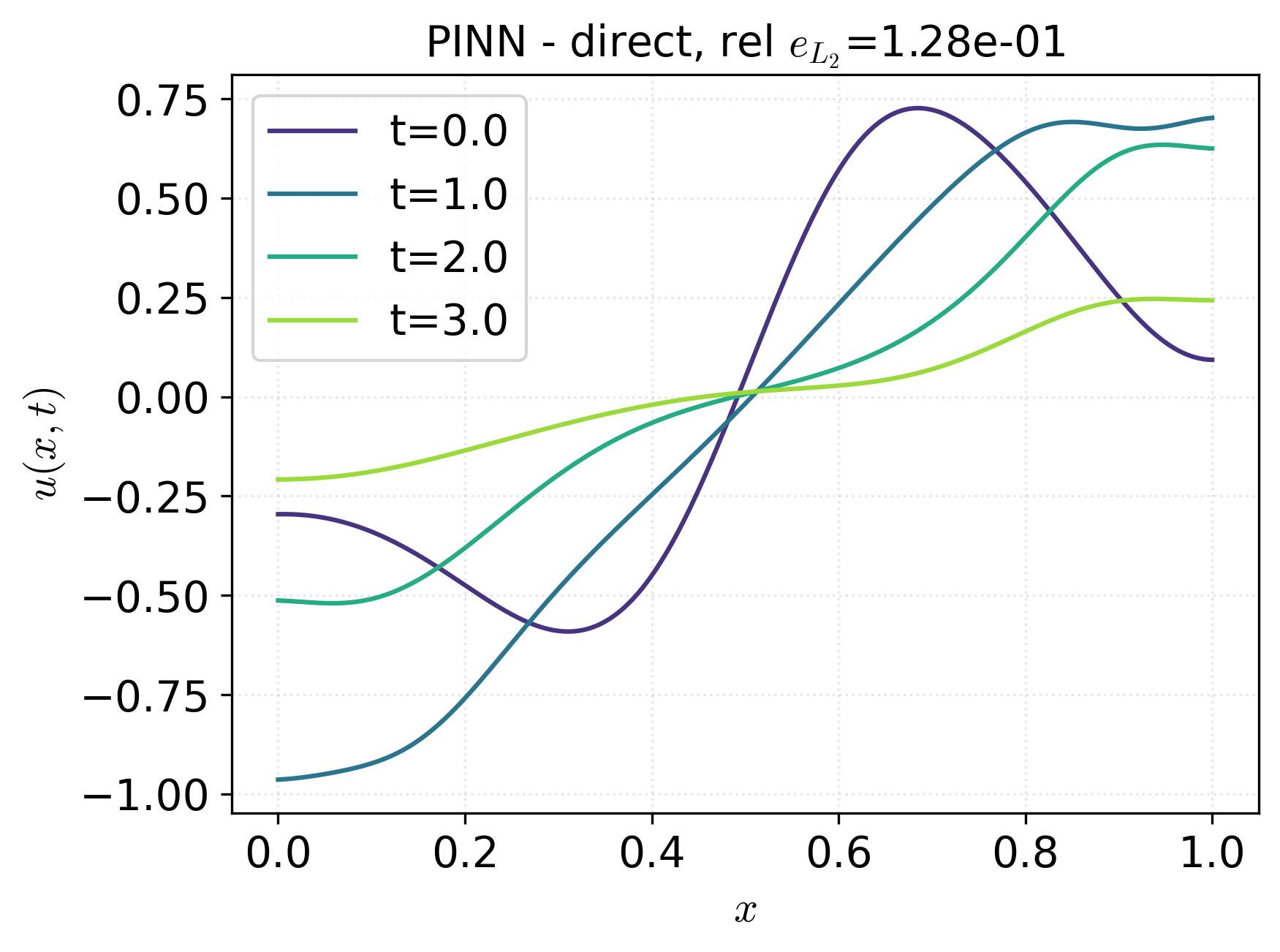}
        \label{fig:ctrl_pinn_direct}
    \end{subfigure}

    \vspace{0.1em}

    \begin{subfigure}[b]{0.7\columnwidth}
        \centering
        \includegraphics[width=\textwidth]{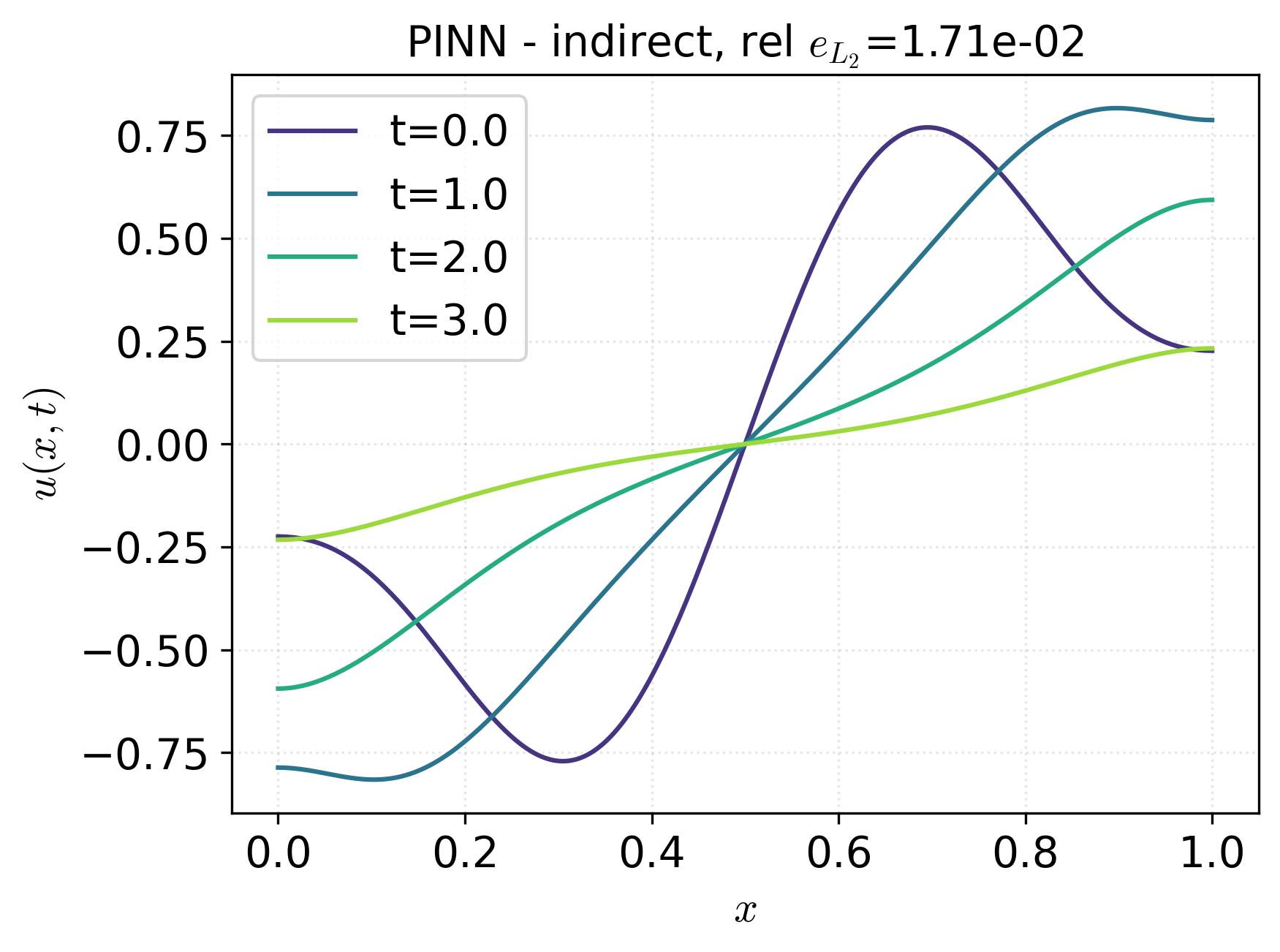}
        \label{fig:ctrl_pinn_indirect}
    \end{subfigure}
    \vspace{0.1em}
    \begin{subfigure}[b]{0.7\columnwidth}
        \centering
        \includegraphics[width=\textwidth]{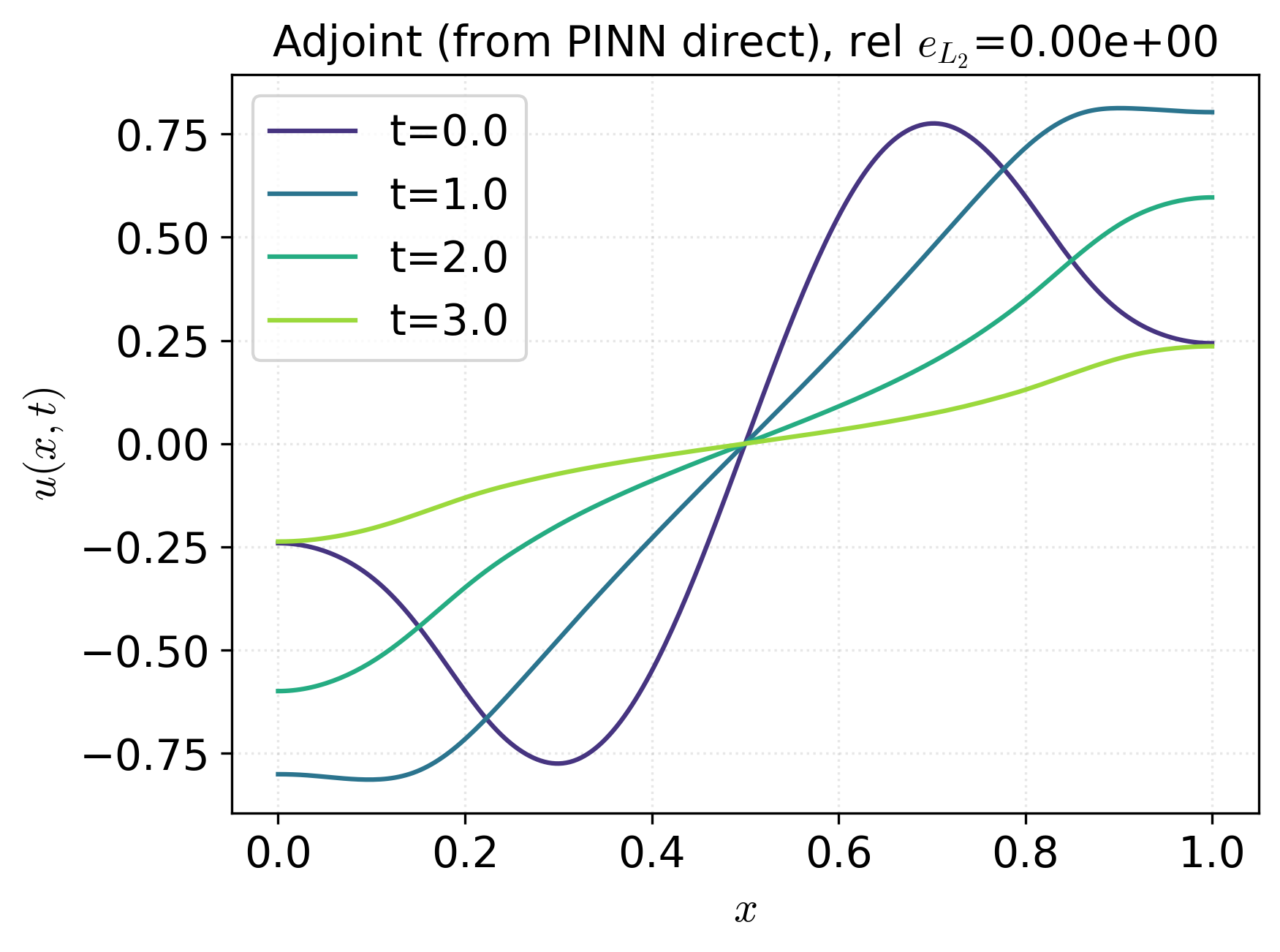}
        \label{fig:ctrl_adjoint_from_pinn_direct}
    \end{subfigure}
    \caption{Control profiles, shown from top to bottom: adjoint starting from scratch, direct PINN, indirect PINN, and adjoint initialized from the direct PINN solution, which is treated as the reference value for this optimal control problem. The relative $L_2$ errors of control are shown in titles. Note that according to the optimality condition, the terminal solution is $y(T)=-\beta_{Q}u(T), \beta_{Q}=1\times10^{-3}$.}
    \label{fig:ctrl_comparison_vertical}
\end{figure}

\begin{figure}[!ht]
    \centering
    \begin{subfigure}[b]{0.7\columnwidth}
        \centering
        \includegraphics[width=\textwidth]{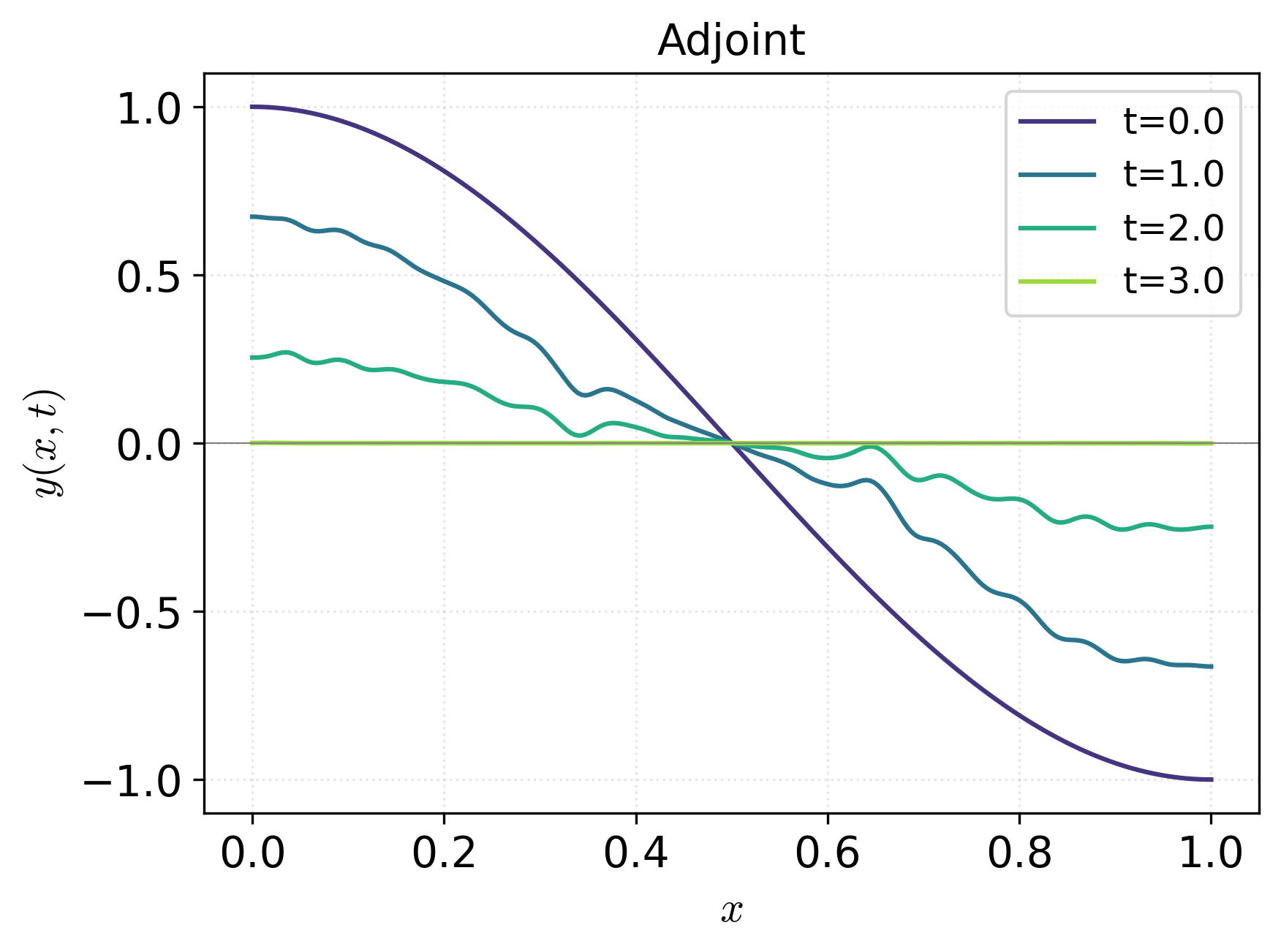}
    \end{subfigure}

    \vspace{0.1em}

    \begin{subfigure}[b]{0.7\columnwidth}
        \centering
        \includegraphics[width=\textwidth]{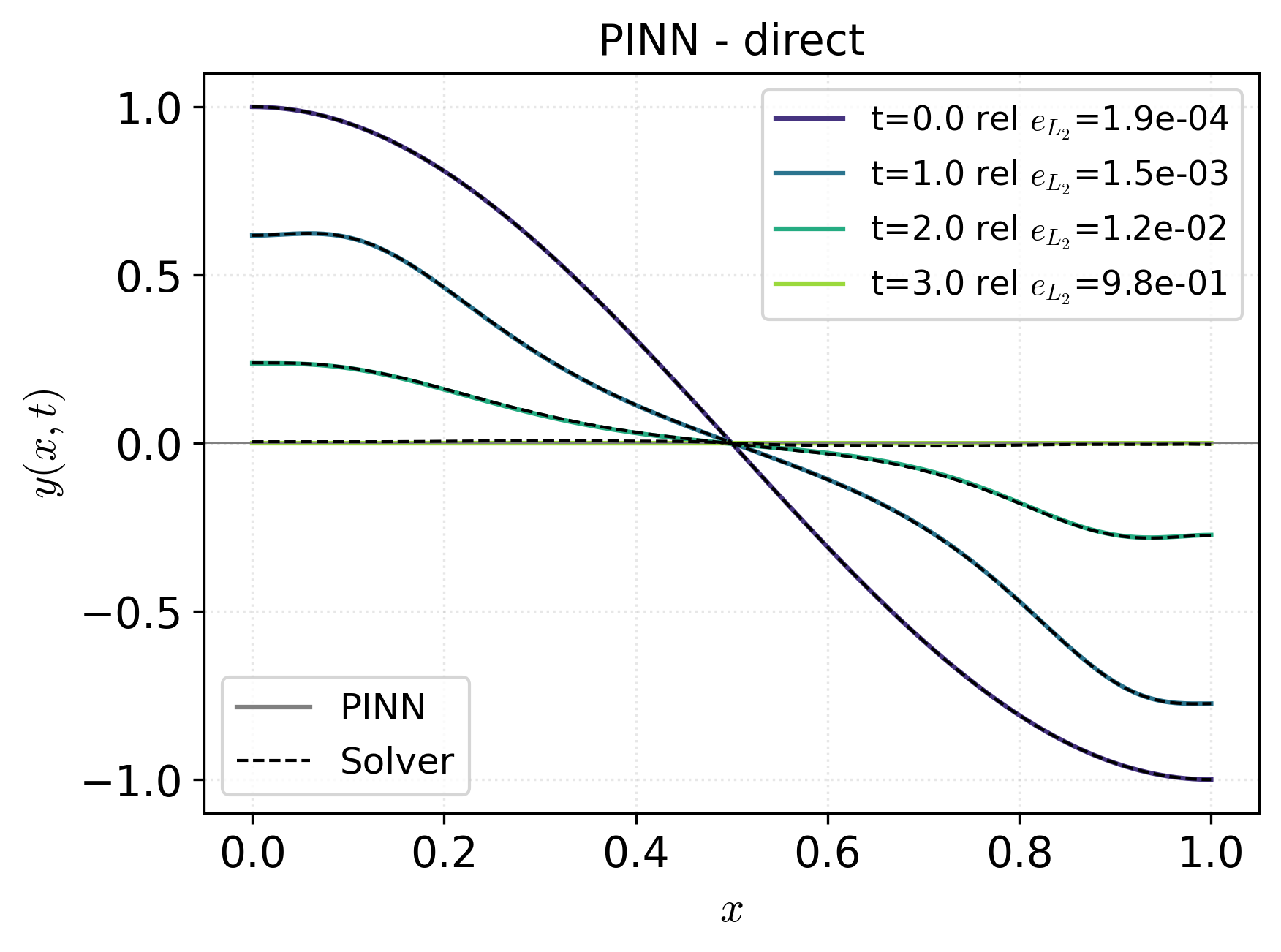}
    \end{subfigure}

    \vspace{0.1em}

    \begin{subfigure}[b]{0.7\columnwidth}
        \centering
        \includegraphics[width=\textwidth]{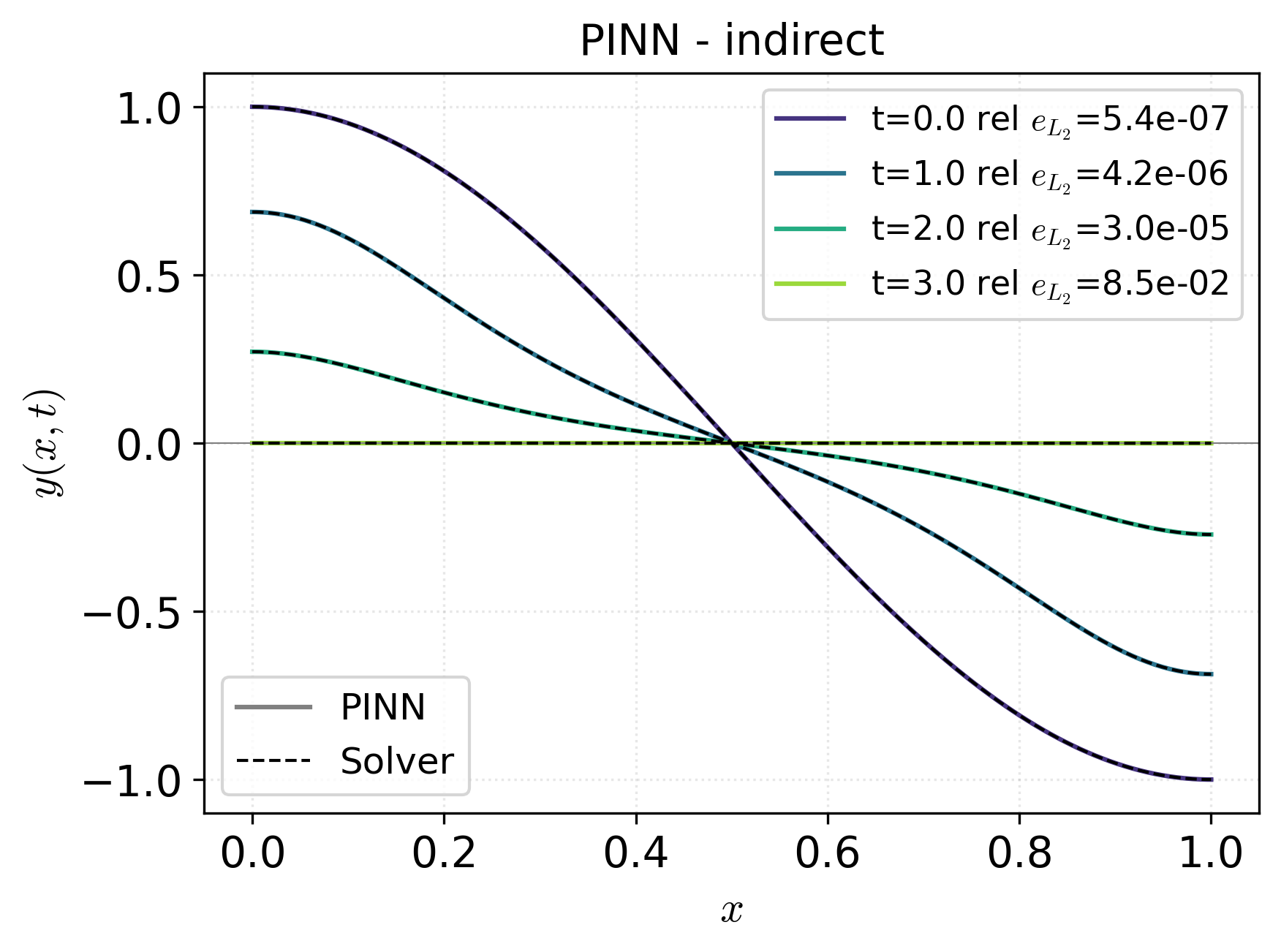}
    \end{subfigure}

    \vspace{0.1em}
    
    \begin{subfigure}[b]{0.7\columnwidth}
        \centering
        \includegraphics[width=\textwidth]{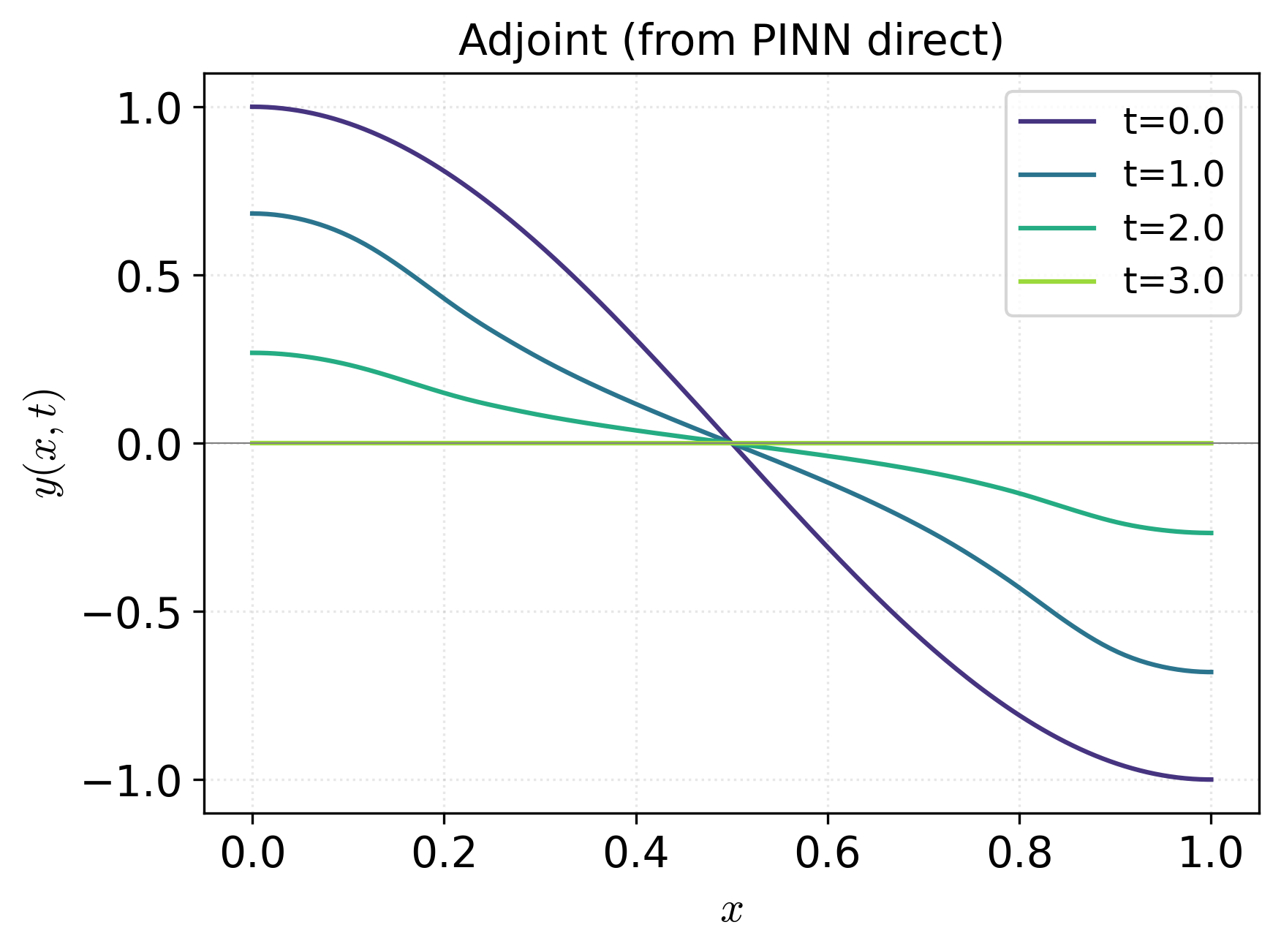}
    \end{subfigure}

    \caption{State trajectories corresponding to the converged controls, shown from top to bottom: adjoint starting from scratch, direct PINN, indirect PINN, and adjoint initialized from the direct PINN solution. For the PINN-based methods, each panel also compares the state predicted directly by the neural network with the state recomputed by the numerical solver using the converged PINN control. The relative $L_2$ error of each time snapshot is marked in the legend.}
    \label{fig:state_comparison_vertical}
\end{figure}

Figures~\ref{fig:losses}--\ref{fig:state_comparison_vertical} highlight the main differences in formulation, convergence, and accuracy among the adjoint, direct PINN, and indirect PINN approaches.

\paragraph{Direct and indirect PINNs solve different optimization problems.}
The key distinction between the two neural formulations is not merely the
training strategy, but the problem being solved. In the direct PINN, the
objective functional is placed directly into the loss together with the PDE,
boundary, and initial-condition residuals. As a result, the optimization is a
penalized surrogate of the original constrained optimal control problem, and
the relative weights of the loss terms effectively modify the problem itself.
In particular, the state equation is no longer enforced exactly, but only in a
soft sense through a weighted residual. Consequently, a small total training
loss does not necessarily imply that the original PDE-constrained problem has
been solved accurately. This effect is clearly visible in
Figs.~\ref{fig:ctrl_comparison_vertical} and
\ref{fig:state_comparison_vertical}: the direct PINN exhibits a relatively
large control error, and its state predicted by the network does not fully
match the state recomputed by the numerical solver using the converged control,
especially at later times. These discrepancies indicate that the objective term
can dominate the weighted loss before the state equation is sufficiently
resolved, and the loss function cannot be reduced to zero, as shown in the loss history of direct PINN in Fig~\ref{fig:losses}.

By contrast, the indirect PINN is trained on the first-order optimality system
itself. Its residual terms correspond to the state equation, adjoint equation,
boundary and initial conditions, and terminal optimality condition, and the optimality condition is satisfied by construction. Therefore,
once these residuals are driven close to zero, the network approaches a true
KKT point of the original optimal control problem, rather than the minimizer of
a penalized surrogate. This explains why the indirect formulation gives
substantially smaller control and state errors and is much closer to the
reference solution, as shown in Figs.~\ref{fig:ctrl_comparison_vertical} and~\ref{fig:state_comparison_vertical}.

\paragraph{Implicit regularization and the role of smooth initialization.}
Another important feature of the PINN solutions is their smoothness. In both
PINN formulations, the control is represented by a smooth low-capacity neural
network with $\tanh$ activations. This parameterization acts as an implicit
regularizer: it suppresses high-frequency oscillations and favors low-complexity,
smooth controls. The effect is evident in
Fig.~\ref{fig:ctrl_comparison_vertical}, where both PINN controls are visibly
smoother than the control obtained by adjoint optimization starting from
scratch. In the discrete adjoint method, the control is represented in a
piecewise-constant form, and once oscillatory components appear, they are not
easily removed during optimization. This is reflected both in the jagged
control profile and in the corresponding fluctuations of the state.

At the same time, the restarted adjoint result reveals that these oscillations
are not intrinsic to the discrete adjoint formulation itself, but are strongly
related to initialization. When the adjoint solver is initialized from the
smooth direct-PINN control, it converges to a much better discrete optimum,
which we treat here as the reference solution. This observation is important:
although the direct PINN is not sufficiently accurate by itself, it provides a
smooth initial guess that places the adjoint iteration in a much better basin
of attraction. In this sense, the PINN acts not only as an optimizer, but also
as a regularized initializer for classical solvers.

\paragraph{Second-order refinement and \texttt{float64} are essential for meaningful PINN convergence.}
The loss histories in Fig.~\ref{fig:losses} show that the indirect PINN
continues to improve long after the total loss has already become very small.
In particular, around iteration $4000$, the total loss has already dropped to
approximately $10^{-9}$, yet the relative control error is still far from its
final value. Continued SSBroyden iterations further reduce the error by nearly
one order of magnitude, from the $10^{-1}$ level to the $10^{-2}$ level. This
late-stage improvement is crucial: it is precisely this regime that separates a
visually plausible PINN solution from a quantitatively meaningful solution for
optimal control.

This behavior highlights two practical requirements. First, second-order or
quasi-Newton refinement is indispensable once the optimization enters the
small-residual regime; first-order training alone would likely stop too early.
Second, double precision is equally important, because the remaining residuals
and parameter updates become too small for robust refinement in single
precision. The present results therefore support a clear practical guideline:
for PDE-constrained PINNs with small or medium-sized networks, a short first-order
warm-up followed by a second-order optimizer in \texttt{float64} is not merely
helpful, but often necessary to obtain a trustworthy solution.

\paragraph{Practical use of direct PINN, indirect PINN, and adjoint refinement.}
From a practical point of view, the three methods serve different purposes.
The direct PINN is the easiest formulation to deploy, since it does not require
deriving the adjoint equation or the full KKT system. It is therefore the most
user-friendly entry point, and it can provide a smooth and physically
reasonable initial guess at relatively low implementation cost. In the present
example, this makes it a natural warm start for either the discrete adjoint
solver or the indirect PINN.

The indirect PINN, however, is the most accurate neural formulation. Because it
targets the optimality system directly, it achieves much better agreement with
the reference control and with the solver-recomputed state. In the current
test, it is also competitive in wall-clock cost and even slightly faster than
the adjoint optimization started from scratch, while giving a smoother and more
accurate solution. More importantly, owing to the implicit regularization of
the neural parameterization, the indirect PINN can outperform the adjoint
method when both are started from naive initial guesses. A practical workflow
suggested by these results is therefore hybrid: use a direct PINN to obtain a
smooth initial control, then refine with either an indirect PINN or a discrete
adjoint solve. This combines ease of use, smoothness, and final accuracy.


\bibliographystyle{IEEEtran}
\bibliography{ref}

\end{document}